\documentclass[a4paper,12pt]{amsart}
\usepackage{amssymb}
\usepackage{ifthen}
 \usepackage[dvips]{graphicx}
\nonstopmode \numberwithin{equation}{section}
\usepackage{amssymb}
\usepackage{ifthen}
\usepackage{graphicx}
\usepackage{amsmath}
\usepackage[T1]{fontenc} 
\usepackage[utf8]{inputenc}
\usepackage[usenames,dvipsnames]{color}
\usepackage{color}
\usepackage[english]{babel}
\usepackage{fancyhdr}
\usepackage{fancybox}
\usepackage{tikz}

\nonstopmode \numberwithin{equation}{section}
\theoremstyle{plain}
\newtheorem{thm}[equation]{Theorem}
\newtheorem{cor}[equation]{Corollary}
\newtheorem{lem}[equation]{Lemma}
\newtheorem{prop}{Proposition}

\newtheorem{conj}{Conjecture}

\theoremstyle{definition}
\newtheorem{defn}{Definition}[section]

\newtheorem{prob}{Problem}
\newtheorem{rem}{Remark}[section]

\newcounter{minutes}
\divide\time by 60
\newcounter{hours}
\multiply\time by 60
\addtocounter{minutes}{-\time}

\newcounter {own}
\def\theown {\thesection       .\arabic{own}}

\newenvironment{pf}[1][]{%
 \vskip 3mm
 \noindent
 \ifthenelse{\equal{#1}{}}%
  {{\slshape Proof. }}%
  {{\slshape #1.} }%
 }%
{\qed\bigskip}

\newcounter{alphabet}

\def\be{\begin{equation}}
\def\ee{\end{equation}}

\newcommand{\bee}{\begin{enumerate}}
\newcommand{\eee}{\end{enumerate}}

\newcommand{\blem}{\begin{lem}}
\newcommand{\elem}{\end{lem}}
\newcommand{\bthm}{\begin{thm}}
\newcommand{\ethm}{\end{thm}}
\newcommand{\bcor}{\begin{cor}}
\newcommand{\ecor}{\end{cor}}
\newcommand{\beg}{\begin{examp}}
\newcommand{\eeg}{\end{examp}}
\newcommand{\begs}{\begin{examples}}
\newcommand{\eegs}{\end{examples}}

\newcommand{\bdefn}{\begin{defn}}
\newcommand{\edefn}{\end{defn}}

\newcommand{\bprob}{\begin{prob}}
\newcommand{\eprob}{\end{prob}}
\newcommand{\bei}{\begin{itemize}}
\newcommand{\eei}{\end{itemize}}

\newcommand{\bcon}{\begin{conj}}
\newcommand{\econ}{\end{conj}}
\newcommand{\bcons}{\begin{conjs}}
\newcommand{\econs}{\end{conjs}}
\newcommand{\bprop}{\begin{prop}}
\newcommand{\eprop}{\end{prop}}
\newcommand{\br}{\begin{rem}}
\newcommand{\er}{\end{rem}}
\newcommand{\brs}{\begin{rems}}
\newcommand{\ers}{\end{rems}}
\newcommand{\bo}{\begin{obser}}
\newcommand{\eo}{\end{obser}}
\newcommand{\bos}{\begin{obsers}}
\newcommand{\eos}{\end{obsers}}
\newcommand{\bpf}{\begin{pf}}
\newcommand{\epf}{\end{pf}}
\newcommand{\ba}{\begin{array}}
\newcommand{\ea}{\end{array}}
\newcommand{\beq}{\begin{eqnarray}}
\newcommand{\beqq}{\begin{eqnarray*}}
\newcommand{\eeq}{\end{eqnarray}}
\newcommand{\eeqq}{\end{eqnarray*}}

\begin{document}

\title{The Bohr Phenomenon for analytic functions on simply connected domains}

\author{Molla Basir Ahamed}
\address{Molla Basir Ahamed,
	School of Basic Science,
	Indian Institute of Technology Bhubaneswar,
	Bhubaneswar-752050, Odisha, India.}
\email{mba15@iitbbs.ac.in}

\author{Vasudevarao Allu}
\address{Vasudevarao Allu,
School of Basic Science,
Indian Institute of Technology Bhubaneswar,
Bhubaneswar-752050, Odisha, India.}
\email{avrao@iitbbs.ac.in}

\author{Himadri Halder}
\address{Himadri Halder,
School of Basic Science,
Indian Institute of Technology Bhubaneswar,
Bhubaneswar-752050, Odisha, India.}
\email{hh11@iitbbs.ac.in}

\subjclass[{AMS} Subject Classification:]{Primary 30C45, 30C50, 30C80}

\def\thefootnote{}
\footnotetext{ {\tiny File:~\jobname.tex,
printed: \number\year-\number\month-\number\day,
          \thehours.\ifnum\theminutes<10{0}\fi\theminutes }
} \makeatletter\def\thefootnote{\@arabic\c@footnote}\makeatother

\begin{abstract}
 In this paper, we investigate the Bohr phenomenon for the class of analytic functions defined on the simply connected domain
 \begin{equation*}
 	\Omega_{\gamma}=\bigg\{z\in\mathbb{C} : \bigg|z+\frac{\gamma}{1-\gamma}\bigg|<\frac{1}{1-\gamma}\bigg\}\;\; \text{for}\;\; 0\leq \gamma<1.
 \end{equation*}
We study improved Bohr radius, Bohr-Rogosinski radius and refined Bohr radius for the class of analytic functions defined in $ \Omega_{\gamma} $, and obtain several sharp results. 
\end{abstract}

\maketitle
\pagestyle{myheadings}
\markboth{Molla Basir Ahamed, Vasudevarao Allu and  Himadri Halder}{On Bohr Phenomenon for simply connected domain}

\section{Introduction and Preliminaries}
\vspace{1cm}
Let $\mathcal{B}(\mathbb{D})$ be the class of analytic functions in unit disk $ \mathbb{D}=\{z\in\mathbb{C} : |z|<1\} $ such that $f(\mathbb{D})\subseteq \overline{\mathbb{D}}$. The classical Bohr theorem for functions $f \in \mathcal{B}(\mathbb{D})$ says that if $f(z)=\sum_{n=0}^{\infty}a_nz^n$, then its associated majorant series $M_{f}(r)$ satisfies the following inequality
\begin{equation}\label{e-1.1}
M_{f}(r):=\sum_{n=0}^{\infty}|a_n|r^n\leq 1\;\; \text{for}\;\; |z|=r\leq\frac{1}{3}
\end{equation}
and the constant $1/3$, called Bohr radius for the class $\mathcal{B}(\mathbb{D})$, cannot be improved. The inequality \eqref{e-1.1} is known as classical Bohr inequality \eqref{e-1.1} for the class $\mathcal{B}(\mathbb{D})$. The Bohr inequality was first obtained by Harald Bohr\;\cite{Bohr-1914} in 1914 with the constant $1/6$. The optimal value $ 1/3 $, which is called the Bohr radius for disk case was later established independently by Weiner, Riesz and Schur. For the proofs we refer to \cite{sidon-1927} and \cite{tomic-1962}. The notion of Bohr inequality has been generalized to several complex variables by finding the multidimensional Bohr radius. We refer the reader to the articles \cite{aizn-2000,aizenberg-2001,boas-1997,Liu-Pon-PAMS-2020}. For more information and intriguing aspects on Bohr phenomenon, we suggest the reader to glance through the articles \cite{Abu-2010}--\cite{Abu - Ali - Pon - 2017}, \cite{aizn-2007}--\cite{alkhaleefah-2019} and \cite{Himadri-Vasu-P1}--\cite{Himadri-Vasu-P4}. Bohr phenomenon for operator valued functions have been extensively studied by Bhowmik and Das (see \cite{Bhowmik-Das--March-2020,Bhowmik-Das-November-2020}).
\vspace{4mm}


The main aim of this article is to study the Bohr inequality for the class of analytic functions that are defined in a general simply connected domain in the complex plain. Let $\Omega$ be a simply connected domain containing $\mathbb{D}$ and $\mathcal{B}(\Omega)$ be the class of analytic functions in $\Omega$ such that $f(\Omega) \subseteq \overline{\mathbb{D}}$. We define the Bohr radius $B=B_{\Omega}$ for the class $\mathcal{B}(\Omega)$ by 
\begin{equation*}
	B:=\sup\bigg\{r\in (0,1) : \sum_{n=0}^{\infty}|a_n|r^n\leq 1\; \text{for all}\; f\in\mathcal{B}(\Omega)\;\; \mbox{with} f(z)=\sum_{n=0}^{\infty}a_nz^n, \;\; z\in\mathbb{D}\bigg\}.
\end{equation*}
In particular, if $ \Omega=\mathbb{D} $, then $ B_{\mathbb{D}}=1/3 $, which is the classical Bohr radius for the class $\mathcal{B}(\mathbb{D})$. Let $\mathbb{D}(a,r):=\{z\in \mathbb{C}: |z-a|<r\}$. Clearly, $\mathbb{D}:=\mathbb{D}(0,1)$. Let $0\leq \gamma <1$. We consider the disk $\Omega _{\gamma}$ defined by 
$$
\Omega_{\gamma}:=\bigg\{z\in\mathbb{C} : \bigg|z+\frac{\gamma}{1-\gamma}\bigg|<\frac{1}{1-\gamma}\bigg\}.
$$
It is easy to see that $\Omega _{\gamma}$ always contains the unit disk $\mathbb{D}$.
In $ 2010 $, the notion of classical Bohr inequality \eqref{e-1.1} has been generalized by Fournier and Ruscheweyh \cite{Four-Rusc-2010} to the class $\mathcal{B}(\Omega_{\gamma})$. More precisely,
\begin{thm}\cite{Four-Rusc-2010} \label{thm-1.2}
	For $ 0\leq \gamma<1 $, let $ f\in\mathcal{B}(\Omega_{\gamma}) $, with $ f(z)=\sum_{n=0}^{\infty}a_nz^n $ in $ \mathbb{D} $. Then,
	\begin{equation*}
		\sum_{n=0}^{\infty}|a_n|r^n\leq 1\;\; \text{for}\;\; r\leq\rho:=\frac{1+\gamma}{3+\gamma}.
	\end{equation*}
Moreover, $ \sum_{n=0}^{\infty}|a_n|\rho^n=1 $ holds for a function $ f(z)=\sum_{n=0}^{\infty}a_nz^n $ in $ \mathcal{B}(\Omega_{\gamma}) $ if, and only if, $ f(z)=c $ with $ |c|=1 $.
\end{thm}
 In this article, we study the Bohr-Rogosinski radius for the class $\mathcal{B}(\Omega_{\gamma})$. In 2017, Kayumov and Ponnusamy \cite{kayumov-2017} introduced Bohr-Rogosinski radius motivated from Rogosinski radius for bounded analytic functions in $\mathbb{D}$. Rogosinski radius is defined as follws: Let $f(z)=\sum_{n=0}^{\infty} a_{n}z^{n}$ be analytic in $\mathbb{D}$ and its corresponding partial sum of $f$ is defined by $S_{N}(z):=\sum_{n=0}^{N-1} a_{n}z^{n}$. Then, for every $N \geq 1$, we have $|\sum_{n=0}^{N-1} a_{n}z^{n}|<1$ in the disk $|z|<1/2$ and the radius $1/2$ is sharp. Motivated by Rogosinski radius, Kayumov and Ponnusamy have considered the Bohr-Rogosinski sum $R_{N}^{f}(z)$ is defined by 
 \begin{equation}
 R_{N}^{f}(z):=|f(z)|+ \sum_{n=N}^{\infty} |a_{n}||z|^{n}.
 \end{equation}
 It is worth to point out that $|S_{N}(z)|=\big|f(z)-\sum_{n=N}^{\infty} a_{n}z^{n}\big| \leq |R_{N}^{f}(z)|$. Thus, it is easy to see that the validity of Bohr-type radius for $R_{N}^{f}(z)$, which is related to the classical Bohr sum (Majorant series) in which $f(0)$ is replaced by $f(z)$, gives Rogosinski radius in the case of bounded analytic functions in $\mathbb{D}$. There has been significant and extensive research carried out on Improved-Bohr inequality and Bohr-Rogosinski radius (see \cite{Alkhaleefah-Kayumov-Ponnusamy-2020,Huang-Liu-Ponnu-2020,Ismagilov-2020,kayumov-2017,kayumov-2018-b,Kay & Pon & AASFM & 2019,kayumov-2018-c,kay & Pon & Sha & MN & 2018,Liu-2020,Ponnusamy-Wirths-2020}).
 \begin{lem}\cite{Ruscheweyh-1985}\label{lem-1.4}
 	Let $ a\in\mathbb{D} $ and $ f\in\mathcal{B}(\mathbb{D}) $ with \begin{equation*}
 	f(z)=\sum_{n=0}^{\infty}\alpha_n(z-a)^n,\;\; |z-a|\leq 1-|a|.
 	\end{equation*}
 	Then, 
 	\begin{equation*}
 	|\alpha_n|\leq (1+|a|)^{n-1}\frac{1-|\alpha_0|^2}{(1-|a|^2)^n},\; n\geq 1.
 	\end{equation*}
 \end{lem}
\noindent Recently, Evdoridis \textit{et al.} \cite{Evd-Ponn-Rasi-2020} obtained the following coefficient bounds for functions defined in $ \Omega_{\gamma} $.
\begin{lem}\cite{Evd-Ponn-Rasi-2020}\label{lem-2.11}
	For $ \gamma\in [0,1) $, let \begin{equation*}
		\Omega_{\gamma}=\bigg\{z\in\mathbb{C} : \bigg|z+\frac{\gamma}{1-\gamma}\bigg|<\frac{1}{1-\gamma}\bigg\},
	\end{equation*}
	and let $ f $ be an analytic function in $ \Omega_{\gamma} $, bounded by $ 1 $, with the series representation $ f(z)=\sum_{n=0}^{\infty}a_nz^n $ in the unit disk $ \mathbb{D} $. Then 
	\begin{equation*}
		|a_n|\leq\frac{1-|a_0|^2}{1+\gamma}\;\; \text{for}\; n\geq 1.
	\end{equation*}
\end{lem}

\section{Main Results}
Before we state an improved version of inequality of Theorem \ref{thm-1.2}, we prove the following lemma. 
\begin{lem}\label{lem-2.2}
	Let $ g : \mathbb{D}\rightarrow \overline{\mathbb{D}} $ be an analytic function, $ m (\geq 2) $ be an integer, and let $ \gamma\in\mathbb{D} $ be such that $ g(z)=\sum_{n=0}^{\infty}\alpha_n(z-\gamma)^n $ for $ |z-\gamma|\leq 1-|\gamma| $. Then
	\begin{equation}\label{e-2.2a}
		\sum_{n=0}^{\infty}\left(|\alpha_n|+\beta|\alpha_n|^m\right)\rho^n\leq 1\;\; \text{for}\;\; \rho\leq \rho_0:=(1-\gamma^2)/(3+\gamma),
	\end{equation} where 
	\begin{equation*}
		\beta=\frac{(1-\gamma)^m(3+\gamma)-(1-\gamma^2)}{8(m-1)}\;\; \text{for}\;\; 0\leq \gamma\leq\gamma_*<1,
	\end{equation*} where $ \gamma_* $ is the smallest root of the equation $ (1-\gamma)^m(3+\gamma)+\gamma^2-1=0. $
\end{lem}
Using Lemma \ref{lem-2.2}, we obtain the following improved version of  Theorem \ref{thm-1.2} for the class $\mathcal{B}(\Omega_{\gamma})$.
\begin{thm}\label{th-3.1}
	For $ 0\leq \gamma<1 $, and integer $ m\;(\geq 2) $, let $ f\in\mathcal{B}(\Omega_{\gamma}) $ with $ f(z)=\sum_{n=0}^{\infty}a_nz^n $ for $ z\in \mathbb{D}, $ then we have 
	\begin{equation*}
		|a_0|+\sum_{n=1}^{\infty}\left(|a_n|+\beta\frac{|a_n|^m}{(1-\gamma)^{(m-1)n}}\right)r^n\leq 1\; \;\text{for}\;\; r\leq r_0=\frac{1+\gamma}{3+\gamma},
	\end{equation*}
where $\beta$ as in Lemma \ref{lem-2.2}. Furthermore, the quantities $ \beta $ and $ (1+\gamma)/(3+\gamma) $ cannot be improved. 
\end{thm}

\begin{figure}[!htb]
	\begin{center}
		\includegraphics[width=0.60\linewidth]{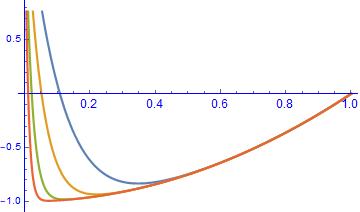}
	\end{center}
	\caption{The roots $ \gamma_*(m) $ of the equation $ (1-\gamma)^m(3+\gamma)+\gamma^2-1=0$.}
\end{figure}
Figure 1 demonstrates values of $ \gamma_* $ in $ [0,1) $ for which $ \beta(\gamma)>0 $ with $ 0\leq \gamma\leq\gamma_*<1 $. The values of $ \gamma_* $ are $ \gamma_*(10)=0.1083,\; \gamma_*(21)=0.0519,\; \gamma_*(50)=0.0219 $ and $ \gamma_*(100)=0.011 $.

\begin{lem}\label{lem-2.7}
	Let $ g : \mathbb{D}\rightarrow\overline{\mathbb{D}} $ be an analytic function, $ \lambda\in [0,512/243] $ and let $ \gamma\in\mathbb{D} $ be such that $ g(z)=\sum_{n=0}^{\infty}\alpha_n(z-\gamma)^n $ for $ |z-\gamma|<1-|\gamma| $. Then 
	\begin{align*}
		\sum_{n=0}^{\infty}|\alpha_n|\rho^n+\left(\frac{8}{9}-\frac{27}{64}\lambda\right)\left(\frac{S^{\gamma}_{\rho}}{\pi}\right)+\lambda\left(\frac{S^{\gamma}_{\rho}}{\pi}\right)^2\leq 1\;\; \text{for}\;\; \rho\leq \rho_0=\frac{1-|\gamma|^2}{3+|\gamma|},
	\end{align*}
	where $ S^{\gamma}_{\rho} $ denotes the area of the image of the disk $ \mathbb{D}(\gamma;r(1-|\gamma|)) $ under the mapping  $ g $.
\end{lem}
\noindent By applying Lemma \ref{lem-2.7}, we obtain the following improved version of Theorem \ref{thm-1.2}.
\begin{thm}\label{th-3.2}
	For $ 0\leq \gamma<1 $ and $ 0\leq \lambda\leq 512/243 $, let $ f\in\mathcal{B}(\Omega_{\gamma}) $ with $ f(z)=\sum_{n=0}^{\infty}a_nz^n $ for $ z\in\mathbb{D}, $ then we have 
	\begin{equation*}
		\sum_{n=0}^{\infty}|a_n|r^n+\left(\frac{8}{9}-\frac{27}{64}\lambda\right)\left(\frac{S_{r(1-\gamma)}}{\pi}\right)+\lambda\left(\frac{S_{r(1-\gamma)}}{\pi}\right)^2\leq 1 \;\; \text{for}\;\; r\leq r_0=\frac{1+\gamma}{3+\gamma}.
	\end{equation*}
Furthermore, the quantities $ 8/9-27\lambda/64 $, $ \lambda $ and $ (1+\gamma)/(3+\gamma) $ cannot be improved. 
\end{thm}
\begin{lem}\label{lem-2.1}
	For $ \gamma\in\mathbb{D} $, let $ g\in\mathcal{B}(\mathbb{D}) $ with $ g(z)=\sum_{n=0}^{\infty}\alpha_n(z-\gamma)^n,\; for $ $ |z-\gamma|\leq 1-|\gamma| $, then 
	\begin{equation*}
		|g(z)|+\sum_{n=N}^{\infty}|\alpha_n|\rho^n\leq 1,\;\; \text{for}\;\; \rho\leq\rho_N,
	\end{equation*}
	where $ \rho_N $ is the root of 
	\begin{equation*}
		2(1+\gamma)\rho^N+ (1+\gamma)(1-\gamma)^{N-1}(\rho-1)(1-\gamma-\rho)=0
	\end{equation*}
	in $(0,1)$.
\end{lem}
\noindent Using Lemma \ref{lem-2.1}, we obtain the following Bohr-Rogosinski radius for the class $\mathcal{B}(\Omega_{\gamma})$.
\begin{thm}\label{th-3.3}
	For $ 0\leq \gamma<1 $ and integer $ N\;(\geq 1) $, let $ f\in\mathcal{B}(\Omega_{\gamma}) $ with $ f(z)=\sum_{n=0}^{\infty}a_nz^n $ for $ z\in\mathbb{D}$. Then, we have 
	\begin{equation*}
		\bigg|f\left(\frac{z-\gamma}{1-\gamma}\right)\bigg|+\sum_{n=N}^{\infty}|a_n|r^n\leq 1 \;\; \text{for}\;\; r\leq r_0=\frac{\rho_N}{1-\gamma},
	\end{equation*}
where $ \rho_N $ is the root of the equation 
\begin{equation}
	2(1+\rho)\rho^N+(1+\gamma)(1-\gamma)^{N-1}(\rho-1)(1-\gamma-\rho)=0.
\end{equation}
	Furthermore, the constant $ {\rho_N}/(1-\gamma) $ cannot be improved. 
\end{thm}

\noindent Using Lemma \ref{lem-2.11}, we establish the following refined Bohr inequality for the class $\mathcal{B}(\Omega_{\gamma})$.
\begin{thm}\label{th-3.5}
For $ 0\leq \gamma<1 $, let $ f\in\mathcal{B}(\Omega_{\gamma}) $ with $ f(z)=\sum_{n=1}^{\infty}a_nz^n $ for $ z\in\mathbb{D}$. Then we have 
\begin{equation*}
	\sum_{n=0}^{\infty}|a_{n+1}|r^n+\left(\frac{1}{1+|a_1|}+\frac{r}{1-r}\right)\sum_{n=2}^{\infty}|a_n|^2r^{2(n-1)}\leq 1 \;\; \text{for}\;\; r\leq r_0=\frac{1+\gamma}{3+\gamma}.
\end{equation*}
The constant $ r_0 $ cannot be improved. 
\end{thm}

\section{Proofs of the Main Results}
\begin{proof}[\bf Proof of the Lemma \ref{lem-2.2}]
	Without loss of generality, we may assume that $ \gamma\in [0,1) $. Using Lemma \ref{lem-1.4}, we obtain
	\begin{align}\label{e-2.3}
		\sum_{n=1}^{\infty}|\alpha_n|\rho^n\leq \frac{1-|\alpha_0|^2}{1+\gamma}\sum_{n=1}^{\infty}\left(\frac{\rho}{1-\gamma}\right)^n=\frac{\left(1-|\alpha_0|^2\right)\rho}{(1+\gamma)(1-\gamma-\rho)}.
	\end{align}
	Further, we have 
	\begin{align}\label{e-2.4}
		\sum_{n=1}^{\infty}|\alpha_n|^m\rho^n\leq\frac{(1-|\alpha_0|^2)^m}{(1+\gamma)^m}\sum_{n=1}^{\infty}\left(\frac{\rho}{(1-\gamma)^m}\right)^n=\frac{\left(1-|\alpha_0|^2\right)^m\rho}{(1+\gamma)((1-\gamma)^m-\rho)}.
	\end{align}
	The series in \eqref{e-2.2a} contains positive terms for $ \beta\geq 0 $. Our aim is to find the smallest value of $ \gamma $ in $ [0,1) $ for which $ \beta\geq 0 $. That is
	\begin{equation*}
		\beta=\frac{(1-\gamma)^m(3+\gamma)-(1-\gamma^2)}{8(m-1)}:=\frac{Q(\gamma)}{8(1-m)}\geq 0,
	\end{equation*} where $ Q(\gamma)=(1-\gamma)^m(3+\gamma)-(1-\gamma^2). $ Clearly, $ \gamma=1 $ is a root of $ Q(\gamma) $. Since $ Q(\gamma) $ is a polynomial such that $ Q(0)=2>0 $ and  $ m\geq 2 $, we have 
	\begin{align*}
		Q\left(\frac{9}{10}\right)=\frac{3.9}{10^m}+\frac{81}{100}-1\leq \frac{84.9}{100}-1=-\frac{15.1}{100}<0.
	\end{align*}
	Therefore, there exists at least one root of $ Q(\gamma) $ in $ (0,1) $. Let $ \gamma_* $ be the smallest root of $ Q(\gamma) $. Then, it is easy to see that $ Q(\gamma)\geq 0 $, and hence $ \beta\geq 0 $ for all $ \gamma\in [0,\gamma_*]. $ A simple computation using \eqref{e-2.3} and \eqref{e-2.4} shows that
	\begin{align} &
		|\alpha_0|+\sum_{n=1}^{\infty}|\alpha_n|\rho^n+\beta\sum_{n=1}^{\infty}|\alpha_n|^m\rho^n\\&\leq\nonumber |\alpha_0|+\frac{\left(1-|\alpha_0|^2\right)\rho}{(1+\gamma)(1-\gamma-\rho)}+\beta\frac{\left(1-|\alpha_0|^2\right)^m\rho}{(1+\gamma)((1-\gamma)^m-\rho)}\\&=\nonumber 1+\Psi_{\gamma}(\rho)\leq 1
	\end{align}
	provided $ \Psi_{\gamma}(\rho)\leq 0, $ where 
	\begin{equation*}
		\Psi_{\gamma}(\rho)=\frac{1-|\alpha_0|^2}{1+\gamma}\left(\frac{\rho}{1-\gamma-\rho}\right)+\beta\left(\frac{1-|\alpha_0|^2}{1+\gamma}\right)^m\left(\frac{\rho}{(1-\gamma)^2-\rho}\right)-(1-|\alpha_0|).
	\end{equation*}
	Since $ (1-\gamma)-\rho>(1-\gamma)^m-\rho $, it is easy to see that $ \Psi_{\gamma}(\rho) $ is an increasing function of $ r $ for $ r<(1-\gamma)^m. $ A simplification shows that 
	\begin{align*}&	\Psi_{\gamma}(\rho)\\&=K\left(1+(1-|\alpha_0|^2)^{m-1}\left(\frac{2\beta \rho}{(1+\gamma)[(1-\gamma)^m-\rho]}+\frac{\phi_{\gamma}(\rho)}{(1-|\alpha_0|^2)^{m-1}}\right)-\frac{2}{1+|\alpha_0|}\right),
	\end{align*}
	where 
	\begin{equation*}
		K=\frac{1-|\alpha_0|^2}{2}\;\; \text{and}\;\;\phi_{\gamma}(\rho)=\frac{2r}{(1+\gamma)(1-\gamma-\rho)}-1.
	\end{equation*}
	\noindent Let $ \rho\leq \rho_0 $ be such that $ \Psi_{\gamma}(\rho)\leq \Psi_{\gamma}(\rho_0) $, and $ \phi_{\gamma}(\rho_0)=0 $. Then, it is easy to see that $ \phi_{\gamma}(\rho_0)=0 $ if, and only if, $ \rho_0=(1-\gamma^2)/(3+\gamma). $
	Therefore, it is enough to prove that $ \Psi_{\gamma}(\rho_0)\leq 0 $ for $ |\alpha_0|\leq 1. $ 
	Let $ \beta=\eta\left((1-\gamma)^m(3+\gamma)-(1-\gamma^2)\right) $, then it is easy to see
	\begin{align*} 
		\Psi_{\gamma}(\rho_0)&=K\left(1+2\eta(1-|\alpha_0|^2)^{m-1}\frac{1-\gamma^2}{(1+\gamma)^m}-\frac{2}{1+|\alpha_0|}\right)\\&:=K G_{\gamma}(|\alpha_0|),
	\end{align*}
	where 
	\begin{equation}\label{e-3.4a}
		G_{\gamma}(x)=1+2\eta A(\gamma)(1-x^2)^{m-1}-\frac{2}{x+1}
	\end{equation}
	and
	\begin{equation*}
	A(\gamma)=\frac{1-\gamma^2}{(1+\gamma)^m}>0\;\;\text{for}\;\;\gamma\in [0,1).
	\end{equation*}
	It now remains to show that $ G_{\gamma}(x)\leq 0 $ for $ \gamma\in [0,1) $ and $ x\in [0,1]. $ Since
	\begin{equation*}
		A^{\prime}(\gamma)=-\frac{2(1+\gamma)\gamma+m(1-\gamma^2)}{(1+\gamma)^{m+1}}\leq 0,\;\;\text{for}\;\;\gamma\in [0,1)
	\end{equation*}
	and $ A(0)=1,\;\; A(1)=0$, it follows that $ A(\gamma) $ is a decreasing function and hence $ A(\gamma)\leq A(0)=1. $ Since $ x\leq 1 $ and $ 0<A(\gamma)\leq 1 $, we have 
	\begin{equation*}
		-A(\gamma)x(1+x)^2(1-x^2)^{m-2}>-4.
	\end{equation*}
	From \eqref{e-3.4a}, we have
	\begin{align*}
		\left(G_{\gamma}(x)\right)^{\prime}&=\frac{2}{(1+x)^2}\left(1-2\eta  A(\gamma)(m-1)x(1+x)^2(1-x^2)^{m-2}\right)\\&\geq \frac{2\left(1-8(m-1)\eta\right)}{(1+x)^2}.
	\end{align*}
	Clearly, $ \left(G_{\gamma}(x)\right)^{\prime}>0 $ for $ x\in (0,1) $ whenever $ \eta\leq 1/(8(m-1)). $ Therefore, $ G_{\gamma}(x) $ is an increasing function on $ [0,1] $ for $ \eta\leq 1/(8(m-1))$.  Equivalently,
	\begin{equation*}
		\beta\leq \frac{(1-\gamma)^m(3+\gamma)-(1-\gamma^2)}{8(m-1)}.
	\end{equation*} In particular, $ G_{\gamma}(x)\leq 0 $ for $ \gamma\in [0,\gamma_*] $ and $ x\in [0,1], $ where $ \gamma_* $ is the smallest root of the equation $ (1-\gamma)^m(3+\gamma)-(1-\gamma^2)=0 $. This completes the proof. 
\end{proof}	

\begin{proof}[\bf Proof of Theorem \ref{th-3.1}]
	For $ 0\leq \gamma<1, $ let 
	\begin{equation*}
		\Omega_{\gamma}=\bigg\{z\in\mathbb{C} : \bigg|z+\frac{\gamma}{1-\gamma}\bigg|<\frac{1}{1-\gamma}\bigg\}
	\end{equation*} and the function $ f : \Omega_{\gamma}\rightarrow\mathbb{D} $ be given by $f(z)=\sum_{n=0}^{\infty} a_{n}z_{n}$. Then the function $ g $ defined by 
\begin{equation*}
	g(z)=f\left(\frac{z-\gamma}{1-\gamma}\right)=\sum_{n=0}^{\infty}\frac{a_n}{(1-\gamma)^n}(z-\gamma)^n\;\;\text{for}\;\; |z-\gamma|<1-\gamma
\end{equation*}
belongs to $ \mathcal{B}(\mathbb{D}) $. Applying Lemma \ref{lem-2.2} to the function $ g $, we obtain
	\begin{equation*}
	|a_0|+\sum_{n=1}^{\infty}\left(\frac{|a_n|}{(1-\gamma)^n}+\beta\left(\frac{|a_n|}{(1-\gamma)^n}\right)^m\right)\rho^n\leq 1\;\; \text{for}\;\; \rho\leq \rho_0=\frac{1-\gamma^2}{3+\gamma}.
	\end{equation*}
 That is 
 \begin{equation*}
 	|a_0|+\sum_{n=1}^{\infty}\left(|a_n|+\beta\frac{|a_n|^m}{(1-\gamma)^{(m-1)n}}\right)\left(\frac{\rho}{1-\gamma}\right)^n\leq 1\;\; \text{for}\;\; \rho\leq \rho_0=\frac{1-\gamma^2}{3+\gamma}
 \end{equation*} 
which is equivalent to 
\begin{equation*}
	|a_0|+\sum_{n=1}^{\infty}\left(|a_n|+\beta\frac{|a_n|^m}{(1-\gamma)^{(m-1)n}}\right)r^n\leq 1\;\; \text{for}\;\; r\leq r_0=\frac{1+\gamma}{3+\gamma},
\end{equation*}
where $ \rho=r(1-\gamma) $ and
\begin{equation*}
	\beta=\frac{(1-\gamma)^m(3+\gamma)-(1-\gamma^2)}{8(m-1)}\;\; \text{for}\;\; 0\leq \gamma\leq\gamma_*<1.
\end{equation*} Here $ \gamma_* $ is the smallest root of the equation $ (1-\gamma)^m(3+\gamma)+\gamma^2-1=0. $ \\[0.6mm] \par
In order to prove the sharpness of the radius, we consider the composition function $ f_{a}=h\circ H $ which maps $ \Omega_{\gamma} $ univalently onto $\mathbb{D}$, where $ H : \Omega_{\gamma}\rightarrow\mathbb{D} $ defined by $ H(z)= (1-\gamma)z+\gamma $ and $ h : \mathbb{D}\rightarrow\mathbb{D} $ with $ h(z)=(a-z)/(1-az), $ for $ a\in (0,1) $. A simple computation shows that 
\begin{equation*}
	f_{a}(z)=\frac{a-\gamma-(1-\gamma)z}{1-a\gamma-a(1-\gamma)}=C_0-\sum_{n=1}^{\infty}C_nz^n\;\;\text{for}\;\; z\in\mathbb{D},
\end{equation*} where $ a\in (0,1) $ and 
\begin{equation*}
	C_0=\frac{a-\gamma}{1-a\gamma}\;\; \text{and}\;\; C_n=\frac{1-a^2}{a(1-a\gamma)}\left(\frac{a(1-\gamma)}{1-a\gamma}\right)^n.
\end{equation*}
A simple computation shows that
\begin{align*} &
	|a_0|+\sum_{n=1}^{\infty}\bigg(|a_n|+\beta\frac{|a_n|^m}{(1-\gamma)^{m-1}n}\bigg)r^n\\&= \frac{a-\gamma}{1-a\gamma}+\sum_{n=1}^{\infty}\bigg(\frac{1-a^2}{a(1-a\gamma)}\left(\frac{a(1-\gamma)}{1-a\gamma}\right)^n+\frac{\beta}{(1-\gamma)^{m-1}n}\frac{(1-a^2)^m}{a^m(1-a\gamma)^m}\left(\frac{a(1-\gamma)}{1-a\gamma}\right)^{mn}\bigg)r^n\\&= \frac{a-\gamma}{1-a\gamma}+\frac{(1+a)(1-a)(1-\gamma)r}{(1-a\gamma)\left((1-a\gamma)-ar(1-\gamma)\right)}+\frac{\beta(1-a)^m(1+a)^mr}{(1-\gamma)^{mn-n-m}(1-a\gamma)^m}\\&=1-(1-a)\Phi_{\gamma}(r),
\end{align*}
where 
\begin{align*} 
	\Phi_{\gamma}(r)&=-\frac{(1+a)(1-\gamma)r}{(1-a\gamma)((1-a\gamma)-ar(1-\gamma))}-\frac{\beta(1-a)^{m-1}(1+a)^mr}{(1-\gamma)^{mn-n-m}(1-a\gamma)^m}-\frac{1}{1-a}\left(\frac{a-\gamma}{1-a\gamma}-1\right)\\&=-\frac{(1+a)(1-\gamma)r}{(1-a\gamma)((1-a\gamma)-ar(1-\gamma))}-\frac{\beta(1-a)^{m-1}(1+a)^mr}{(1-\gamma)^{mn-n-m}(1-a\gamma)^m}-\frac{1+\gamma}{1-a\gamma}.
\end{align*}
Therefore, $ \Phi_{\gamma}(r) $ is strictly decreasing function of $ r$ in $(0,1)$. Hence, for $ r>r_0=(1+\gamma)/(3+\gamma) $, we have $ \Phi_{\gamma}(r)<\Phi_{\gamma}(r_0) $. A simple computation shows that 
\begin{align*}
	\lim_{a\rightarrow 1}\Phi_{\gamma}(r_0)=-\frac{2r_0}{(1-\gamma)(1-r_0)}+\frac{1+\gamma}{1-\gamma}=0.
\end{align*}
Thus, $ \Phi_{\gamma}(r)< 0 $ for $ r>r_0 $. Hence, $ 1-(1-a)\Phi_{\gamma}(r)>1 $ for $ r>r_0, $ which shows that $ r_0 $ is the best possible. This completes the proof.
\end{proof}	

\begin{pf}[\bf Proof of Lemma \ref{lem-2.7}]
	Without loss of generality, we assume that $ \gamma\in [0,1) $. Also let $ z\in\mathbb{D}_{\gamma}:=\mathbb{D}(\gamma;1-\gamma) $ if, and only if, $ w=(z-\gamma)/(1-\gamma)\in\mathbb{D}. $
	Then we have 
	\begin{align*}
		g(z)=\sum_{n=0}^{\infty}\alpha_n(1-\gamma)^n\phi^n(z)=\sum_{n=0}^{\infty}b_n\phi^n(z):=G(\phi(z))
	\end{align*}
	for $ z\in \mathbb{D}_{\gamma}, $ where $ b_n= \alpha_{n} (1-\gamma)^n$. A simple computation shows that
	\begin{align}\label{e-2.8}
		\frac{S^{\gamma}_{\rho}}{\pi}=\frac{1}{\pi}\text{Area}\bigg(G(\mathbb{D}(0,\rho))\bigg)\leq (1-|b_0|^2)^2\frac{\rho^2}{(1-\rho^2)^2}=(1-|\alpha_0|^2)^2\frac{\rho^2}{(1-\rho^2)^2}.
	\end{align}
	Therefore,
	\begin{align}\label{e-2.9}
		\sum_{n=1}^{\infty}|\alpha_n|\rho^n\leq\frac{1-|\alpha_0|^2}{1+\gamma}\sum_{n=1}^{\infty}\left(\frac{\rho}{1-\gamma}\right)^n=\frac{1-|\alpha_0|^2}{1+\gamma}\frac{\rho}{1-\gamma-\rho}.
	\end{align}
	In view of \eqref{e-2.8} and \eqref{e-2.9}, we obtain 
	\begin{align*} &
		|\alpha_0|+	\sum_{n=1}^{\infty}|\alpha_n|\rho^n+k\left(\frac{S^{\gamma}_{\rho}}{\pi}\right)+\lambda\left(\frac{S^{\gamma}_{\rho}}{\pi}\right)^2\\&= |\alpha_0|+\frac{(1-|\alpha_0|^2)\rho}{(1+\gamma)(1-\gamma-\rho)}+k\frac{(1-|\alpha_0|^2)^2\rho^2}{(1-\rho^2)^2}+\lambda\frac{(1-|\alpha_0|^2)^4\rho^4}{(1-\rho^2)^4}\\&= 1+\Psi_1^{\gamma}(\rho),
	\end{align*}
	where 
	\begin{align*}
		\Psi_1^{\gamma}(\rho)=\frac{(1-|\alpha_0|^2)\rho}{(1+\gamma)(1-\gamma-\rho)}+k\frac{(1-|\alpha_0|^2)^2\rho^2}{(1-\rho^2)^2}+\lambda\frac{(1-|\alpha_0|^2)^4\rho^4}{(1-\rho^2)^4}-(1-|\alpha_0|)
	\end{align*}
	which can be written as 
	\begin{align*} 
		\Psi_1^{\gamma}(\rho)&=\frac{1-|\alpha_0|^2}{2}\bigg(1+2\lambda(1-|\alpha_0|^2)^3 \left(\frac{\rho^4}{(1-\rho^2)^4}+\frac{k}{\lambda}\frac{\rho^2}{(1-\rho^2)^2(1-|\alpha_0|^2)^2}\right)\\&\;\;\;\;+\frac{1}{2\lambda(1-|\alpha_0|^2)^3}\left(\frac{2\rho}{(1+\gamma)(1-\gamma-\rho)}-1\right)-\frac{2}{1+|\alpha_0|}\bigg).
	\end{align*}
	Let $ \rho\leq \rho_0 $. Then, it is easy to see that $\Psi_1^{\gamma}(\rho)$ is an increasing function and hence $ 	\Psi_1^{\gamma}(\rho)\leq 	\Psi_1^{\gamma}(\rho_0) $, where 
	\begin{align*}
		\frac{2\rho_0}{(1+\gamma)(1-\gamma-\rho_0)}=1,\;\;\text{\it i.e.,}\;\; \rho_0=\frac{1-\gamma^2}{3+\gamma}.
	\end{align*}
	A simple computation shows that 
	\begin{align*}
		\Psi_1^{\gamma}(\rho_0)&=\frac{1-|\alpha_0|^2}{2}\left(1+2\lambda(1-|\alpha_0|^2)^3A^4(\gamma)+2k(1-|\alpha_0|^2)A^2(\gamma)-\frac{2}{1+|\alpha_0|}\right)\\&=\frac{1-|\alpha_0|^2}{2}J(|\alpha_0|),
	\end{align*}
	where 
	\begin{align*}
		J(x)&=1+2\lambda(1-x^2)^3A^4(\gamma)+2k(1-x^2)A^2(\gamma)-\frac{2}{1+x}\;\; \mbox{for}\;\ x\in [0,1]\\\text{and}\;\; A(\gamma)&=\frac{(3+\gamma)(1-\gamma^2)}{(3+\gamma)^2-(1-\gamma^2)^2}.
	\end{align*}
	It is enough to show that $ J(x)\leq 0 $ for $ x\in [0,1] $ and $ \gamma\in [0,1) $ so that $ 	\Psi_1^{\gamma}(\rho_0)\leq 0 $. We note that $ A(\gamma)>0 $ for $ \gamma\in [0,1) $. Further,
	\begin{align*}
		J(0)=2\lambda A^4(\gamma)+2kA^2(\gamma)-1, \;\; \text{and}\;\; \lim_{x\rightarrow 1^{-}}J(x)=0.
	\end{align*}
	It can be seen that $ A(\gamma)=(f_1\circ f_2)(\gamma) $, where $ f_1(\rho)=\rho/(1-\rho^2) $ and $ f_2(\gamma)=(1-\gamma^2)/(3+\gamma). $
	Since $ A^{\prime}(\gamma)=f^{\prime}_1(f_2(\gamma)f^{\prime}_2(\gamma)), $ where 
	\begin{equation}
		f^{\prime}_2(\gamma)=-\left(\frac{\gamma^2+6\gamma+1}{(3+\gamma)^2}\right)<0
	\end{equation}
	which implies that $ f_1(\rho) $ is an increasing function of $ \rho $ in $ (0,1) $, and $ f_2 $ is a decreasing function of $ \gamma $ in $ [0,1) $. Hence, it follows that $ A(\gamma) $ is a decreasing function of $ \gamma $  in $ [0,1) $, with $ A(0)=3/8 $ and $ A(1)=0 $. It can be seen that $ A^2(\gamma) $ and $ A^4(\gamma) $ are decreasing functions on $ [0,1) $. Therefore, we have 
	\begin{align*}
		A^2(\gamma)\leq A^2(0)=\frac{9}{64}\;\; \text{and}\;\; A^4(\gamma)\leq A^4(0)=\frac{81}{4096}.
	\end{align*}
	Since $ x\in [0,1] $, we have
	\begin{align*}
		x(1+x)^2A^2(\gamma)\leq \frac{9}{16}\;\;\text{and}\;\; x(1+x)^2(1-x^2)^2A^4(\gamma)\leq \frac{81}{1024}.
	\end{align*}
	As a consequence, we obtain
	\begin{align*}
		J^{\prime}(x)&=\frac{2}{(1+x)^2}\bigg(1-2kx(1+x)^2A^2(\gamma)-6\lambda x(1+x)^2(1-x^2)^2A^4(\gamma)\bigg)\\&\geq  \frac{2}{(1+x)^2}\bigg(1-\left(\frac{9k}{8}+\frac{243\lambda}{512}\right)\bigg)\\&\geq 0,\;\;\;\;\; \text{if}\;\; k+27\lambda/64\leq 8/9.
	\end{align*}
	Therefore, $ J(x) $ is an increasing function in $ [0,1] $ for $ k+27\lambda/64\leq 8/9 $. Hence, $ J(x)\leq 0 $ for all $ x\in [0,1] $ and $ \gamma\in [0,1) $. This completes the proof.
\end{pf}

\begin{pf}[\bf Proof of Theorem \ref{th-3.2}]
	Let $f \in \mathcal{B}(\Omega_{\gamma})$ and $ g(z)=f((z-\gamma)/(1-\gamma)) $. Then, it is easy to see that $ g\in\mathcal{B}(\mathbb{D})$ and 
	\begin{align*}
		g(z)=\sum_{n=0}^{\infty}\frac{a_n}{(1-\gamma)^n}(z-\gamma)^n.
	\end{align*}
Using Lemma \ref{lem-2.7}, we obtain 
\begin{align*}
	\sum_{n=0}^{\infty}\frac{|a_n|}{(1-\gamma)^n}\rho^n+\left(\frac{8}{9}-\frac{27}{64}\lambda\right)\left(\frac{S^{\gamma}_{\rho}}{\pi}\right)+\lambda\left(\frac{S^{\gamma}_{\rho}}{\pi}\right)^2\leq 1\;\; \text{for}\;\; \rho\leq\frac{1-\gamma^2}{3+\gamma}
\end{align*} which is equivalent to 
\begin{align}\label{e-4.1}
	\sum_{n=0}^{\infty}{|a_n|}\left(\frac{\rho}{(1-\gamma)}\right)^n+\left(\frac{8}{9}-\frac{27}{64}\lambda\right)\left(\frac{S^{\gamma}_{\rho}}{\pi}\right)+\lambda\left(\frac{S^{\gamma}_{\rho}}{\pi}\right)^2\leq 1\;\; \text{for}\;\; \rho\leq\frac{1-\gamma^2}{3+\gamma}.
\end{align}
Set $ \rho=r(1-\gamma) $, then in view of \eqref{e-4.1}, we obtain 
\begin{align}\label{e-3.8e}
	\sum_{n=0}^{\infty}{|a_n|}r^n+\left(\frac{8}{9}-\frac{27}{64}\lambda\right)\left(\frac{S^{\gamma}_{r(1-\gamma)}}{\pi}\right)+\lambda\left(\frac{S^{\gamma}_{r(1-\gamma)}}{\pi}\right)^2\leq 1\;\; \text{for}\;\; r\leq\frac{1+\gamma}{3+\gamma}.
\end{align} 
To show the sharpness of the result, we consider the following function
\begin{align*}
	f_a(z)=\frac{a-\gamma-(1-\gamma)z}{1-a\gamma-a(1-\gamma)}\;\;\text{for}\;\; z\in \Omega_{\gamma}\;\; \text{and}\;\; a\in (0,1).
\end{align*}
Define $ \phi : \mathbb{D}\rightarrow\mathbb{D} $ by $ \phi(z)=(a-z)/(1-az) $ and $ H :\Omega_{\gamma}\rightarrow\mathbb{D} $ by $ H(z)=(1-\gamma)z+\gamma $. 
Then, the function $ f_a=\phi\circ H $ maps $ \Omega_{\gamma} $, univalently onto $ \mathbb{D}. $ A simple computation shows that
\begin{equation*}
	f_a(z)=\frac{a-\gamma-(1-\gamma)z}{1-a\gamma-a(1-\gamma)}=C_0-\sum_{n=1}^{\infty}C_nz^n\;\; \mbox{for}\;\; z\in\mathbb{D},
\end{equation*} where $ a\in (0,1) $ and 
\begin{equation*}
	C_0=\frac{a-\gamma}{1-a\gamma}\;\; \text{and}\;\; C_n=\frac{1-a^2}{a(1-a\gamma)}\left(\frac{a(1-\gamma)}{1-a\gamma}\right)^n.
\end{equation*}
A simple computation using \eqref{e-3.8e} shows that
\begin{align*} &
	\sum_{n=0}^{\infty}{|a_n|}r^n+\left(\frac{8}{9}-\frac{27}{64}\lambda\right)\left(\frac{S^{\gamma}_{r(1-\gamma)}}{\pi}\right)+\lambda\left(\frac{S^{\gamma}_{r(1-\gamma)}}{\pi}\right)^2\\&=\frac{a-\gamma}{1-a\gamma}+\left(\frac{1-a^2}{1-a\gamma}\right)\frac{(1-\gamma)r}{1-a\gamma-ar(1-\gamma)}+\left(\frac{8}{9}-\frac{27}{64}\lambda\right)\frac{r^2(1-a^2)^2(1-\gamma)^4}{((1-a\gamma)^2-a^2r^2(1-\gamma)^4)^2}\\& \;\;\;\;\;\;+\lambda\frac{r^4(1-a^2)^4(1-\gamma)^8}{((1-a\gamma)^2-a^2r^2(1-\gamma)^4)^4} \\& :=1-(1-a)\Phi^{\gamma}_1(r),
\end{align*}
where 
\begin{align*} 
	\Phi^{\gamma}_1(r)&=-\frac{(1+a)(1-\gamma)r}{(1-a\gamma-ar(1-\gamma))(1-a\gamma)}-\left(\frac{8}{9}-\frac{27}{64}\lambda\right)\frac{r^2(1-a)(1+a)^2(1-\gamma)^4}{((1-a\gamma)^2-a^2r^2(1-\gamma)^4)^2}\\&\quad\quad -\lambda\frac{r^4(1-a)^3(1+a)^4(1-\gamma)^8}{((1-a\gamma)^2-a^2r^2(1-\gamma)^4)^4}-\frac{1}{1-a}\left(\frac{a-\gamma}{1+a\gamma}-1\right).
\end{align*}
It is easy to see that $ \Phi^{\gamma}_1(r) $ is strictly decreasing function of $ r$ in $(0,1) $. Therefore,  for $ r>r_0=(1+\gamma)/(3+\gamma) $, we have $ \Phi^{\gamma}_1(r)<\Phi^{\gamma}_1(r_0) $. An elementary calculation shows that 
\begin{align*}
	\lim_{a\rightarrow 1}\Phi^{\gamma}_1(r_0)=-\frac{2r_0}{(1-\gamma)(1-r_0)}+\frac{1+\gamma}{1-\gamma}=0.
\end{align*}
Therefore, $ \Phi^{\gamma}_1(r)< 0 $ for $ r>r_0 $. Hence, $ 1-(1-a)\Phi_{\gamma}(r)>1 $ for $ r>r_0, $ which shows that $ r_0 $ is the best possible. 
\end{pf}
\begin{pf}[\bf Proof of Lemma \ref{lem-2.1}]
	Let $ g\in\mathcal{B}(\mathbb{D}) $. Then by the 
	Schwarz-Pick lemma, for any $ g\in\mathcal{B}(\mathbb{D}) $, we have
	\begin{equation}\label{e-3.8a}
		|g(z)|\leq \frac{\rho+|g(0)|}{1+\rho|g(0)|}\quad \mbox{for}\quad z\in\mathbb{D}.
	\end{equation}
	For functions $ g\in\mathcal{B}(\mathbb{D}) $, from Lemma \ref{lem-1.4}, we have
	\begin{equation}\label{e-3.8b}
		|\alpha_n|\leq (1+|\gamma|)^{n-1}\frac{1-|\alpha_0|^2}{(1-|\gamma|^2)^n}\quad\mbox{for}\quad n\geq 1.
	\end{equation}
	A simple computation using \eqref{e-3.8b} gives
	
	\begin{align}\label{e-3.8c}
		\sum_{n=N}^{\infty}|\alpha_n|\rho^n\leq \frac{1-|\alpha_0|^2}{1+\gamma}\sum_{n=N}^{\infty}\left(\frac{\rho}{1-\gamma}\right)^n=\frac{(1-|\alpha_0|^2)}{(1+\gamma)(1-\gamma)^{N-1}}\left(\frac{\rho^N}{1-\gamma-\rho}\right).
	\end{align}
	From \eqref{e-3.8a} and \eqref{e-3.8c} we obtain
	\begin{align*}
		|g(z)|+	\sum_{n=N}^{\infty}|\alpha_n|\rho^n&\leq \frac{\rho+|g(0)|}{1+\rho|g(0)|}+\frac{(1-|\alpha_0|^2)}{(1+\gamma)(1-\gamma)^{N-1}}\left(\frac{\rho^N}{1-\gamma-\rho}\right)\\&=1+\frac{\Phi^{\gamma}_N(\rho)}{(1+\rho|\alpha_0|)(1+\gamma)(1-\gamma)^{N-1}(1-\gamma-\rho)},
	\end{align*}
	where 
	\begin{align*}
		\Phi^{\gamma}_N(\rho)&=(\rho+|\alpha_0|)A(\gamma)(1-\gamma-\rho)+(1+|\alpha_0|)(1-|\alpha_0|)(1+\rho|\alpha_0|)\rho^N\\&\quad\quad-(1+\rho|\alpha_0|)A(\gamma)(1-\gamma-\rho)\\&=(1-|\alpha_0|)\bigg((1+|\alpha_0|)(1+\rho|\alpha_0|)\rho^N+A(\gamma)(\rho-1)(1-\gamma-\rho)\bigg)\\&\leq (1-|\alpha_0|)\bigg(2(1+\gamma)\rho^N+A(\gamma)(\rho-1)(1-\gamma-\rho)\bigg), 
	\end{align*}
	where  $ A(\gamma)=(1+\gamma)(1-\gamma)^{N-1} $ and $ |\alpha_0|\leq 1 $.
	\noindent An observation shows that $ \Phi^{\gamma}_N(\rho)\leq 0 $ if $ 2(1+\gamma)\rho^N+A(\gamma)(\rho-1)(1-\gamma-\rho)\leq 0 $, and this holds for $ \rho\leq \rho_N $, where $ \rho_N $ is the root of 
	\begin{equation*}
		F_N(\gamma,\rho)=2(1+\gamma)\rho^N+A(\gamma)(\rho-1)(1-\gamma-\rho)=0.
	\end{equation*}
	The existence of the root $ \rho_N $ in $ (0,1) $ follows from the fact that $ F_N(\gamma,\rho) $ is continuous and $ F_N(\gamma,0)F_N(\gamma,1)<0. $ 
\end{pf}

\begin{pf}[\bf Proof of Theorem \ref{th-3.3}]
	For $ 0\leq \gamma<1 $, let $ f\in \mathcal{B}(\Omega_{\gamma}) $ such that $ f(z)=\sum_{n=0}^{\infty}a_nz^n $ for $ z\in\mathbb{D}. $ Then, it is easy to see that  
	\begin{equation*}
		g(z)=f\left(\frac{z-\gamma}{1-\gamma}\right)\in\mathcal{B}(\mathbb{D})\;\; \text{ for}\;\; |z-\gamma|<1-|\gamma|.
	\end{equation*}
Further,
\begin{equation*}
	g(z)=f\left(\frac{z-\gamma}{1-\gamma}\right)=\sum_{n=0}^{\infty}\frac{a_n}{(1-\gamma)^n}(z-\gamma)^n.
\end{equation*}
An application of Lemma \ref{lem-2.1} shows that
\begin{align} \label{e-3.8}	\bigg|f\left(\frac{z-\gamma}{1-\gamma}\right)\bigg|+\sum_{n=N}^{\infty}\frac{|a_n|}{(1-\gamma)^n}\rho^n&\leq 1\;\; \text{for}\;\; \rho\leq\rho_N.
\end{align} 
Since $|z-\gamma|<1-\gamma$, we set  $z-\gamma=w(1-\gamma)$ for some $w \in \mathbb{D}$ and $\rho = r(1-\gamma)$. Then, from \eqref{e-3.8}, we obtain 
\begin{align*}
|f(w)|+ \sum_{n=N}^{\infty} |a_{n}|r^n \leq 1 \;\; \mbox{for} \;\; r\leq \frac{\rho_{N}}{1-\gamma},
\end{align*} 
where $ \rho_N $ as in Lemma \ref{lem-2.1}. That is, $ \rho_N $ is the smallest root of the equation $ 2(1+\rho)\rho^N+A(\gamma)(\rho-1)(1-\gamma-\rho)=0. $
\\ In order to show the sharpness of the result, we consider the following function $ f_a $ defined by 
\begin{equation*}
	f_a(z)=\frac{1-\gamma-(1-\gamma)z}{(1-a\gamma)-(1-\gamma)z}=B_0-\sum_{n=1}^{\infty}B_nz^n\quad\mbox{for}\quad z\in\mathbb{D}.
\end{equation*} 
For $ \gamma\in [0,1) $, $ a>\gamma $ and $ \rho=r(1-\gamma) $, we obtain
\begin{align}\label{e-3.13}
	M:&=|f_a(-\rho)|+\sum_{n=N}^{\infty}|a_n|\rho^n\\&\nonumber=\frac{(a-\gamma)+(1-\gamma)\rho}{(1-a\gamma)+a(1-\gamma)\rho}+\frac{1-a^2}{a(1-a\gamma)}\left(\frac{a(1-\gamma)}{1-a\gamma}\right)^N\rho^N\left(\frac{1-a\gamma}{(1-a\gamma)-a(1-\gamma)\rho}\right)\\&\nonumber=\frac{(a-\gamma)+(1-\gamma)\rho}{(1-a\gamma)+a(1-\gamma)\rho}+\frac{(1-a^2)B^N\rho^N}{a\left((1-a\gamma)-a(1-\gamma)\rho\right)},\;\;\;\; \text{where}\; B=\frac{a(1-\gamma)}{1-a\gamma}\\&\nonumber=\frac{\left((a-\gamma)+(1-\gamma)\rho\right)\left((1-a\gamma)-a(1-\gamma)\rho\right)+d\rho^{N}(1-a^2)\left((1-a\gamma)+a(1-\gamma)\rho\right)}{\left((1-a\gamma)+a(1-\gamma)\rho\right)\left((1-a\gamma)-a(1-\gamma)\rho\right)}
\end{align}
From \eqref{e-3.13}, it is easy to see that $ M>1 $ if $ V(\rho)>0 $, where
\begin{align*} 	V(\rho)&=\left((a-\gamma)+(1-\gamma)\rho\right)\left((1-a\gamma)-a(1-\gamma)\rho\right)+d\rho^{N}(1-a^2)\left((1-a\gamma)+a(1-\gamma)\rho\right)\\&\quad\quad-\left((1-a\gamma)+a(1-\gamma)\rho\right)\left((1-a\gamma)-a(1-\gamma)\rho\right)\\&= (1-a)\bigg((1+a)\left((1-a\gamma)+a(1-\gamma)\rho\right)d\rho^N\\&\quad\quad+\bigg((1-a\gamma)-a(1-\gamma)\rho\bigg)\bigg(\rho(1-\gamma)-(1+\gamma)\bigg)\bigg).
\end{align*}
Note that $ V(\rho)>0 $ if 
\begin{align} \label{e-3.14}
	W(\rho)&:=(1+a)\left((1-a\gamma)+a(1-\gamma)\rho\right)d\rho^N\\ \nonumber&\;\;\;+\bigg((1-a\gamma)-a(1-\gamma)\rho\bigg)\bigg(\rho(1-\gamma)-(1+\gamma)\bigg)>0.
\end{align}
Therefore, $M \leq 1$ for all $a\in [0,1)$, only in the case when $\rho \leq \rho_{N}$. Finally, allowing $a \rightarrow 1$, from the inequality \eqref{e-3.14}, it can be seen that $M>1$ if $\rho> \rho_{N}$. Thus, $M>1$ if $r>\rho_{N}/(1-\gamma)$. This proves the sharpness.

\end{pf}
\begin{pf}[\bf Proof of Theorem \ref{th-3.5}]
	Let $ f\in\mathcal{B}(\Omega_{\gamma}) $ be given by $ f(z)=\sum_{n=1}^{\infty}a_nz^n $ for $ z\in\mathbb{D} $. Then, $ f $ can be expressed as $ f(z)=zh(z) $, where $ h\in B(\Omega_{\gamma}) $ with $ h(z)=\sum_{n=0}^{\infty}b_nz^n $ and $ b_n=a_{n+1}. $ Let  $ |b_0|=|a_1|=a $, and $ h_0(z)=g(z)-b_0 $. Using Lemma \ref{lem-2.11}, we obtain
\begin{align} \label{e-4.2}&
\sum_{n=0}^{\infty}|b_n|r^n+\left(\frac{1}{1+|b_0|}+\frac{r}{1+r}\right)\sum_{n=1}^{\infty}|b_n|^2r^{2n}\\&\leq a+\frac{1-a^2}{1+\gamma}\frac{r}{1-r}+\left(\frac{1}{1+a}+\frac{r}{1-r}\right)\left(\frac{1-a^2}{1+\gamma}\right)^2\frac{r^2}{1-r^2}\nonumber.
\end{align}
That is,
 \begin{align}\label{e-4.3}
	\sum_{n=0}^{\infty}|b_n|r^n&\leq a+\frac{1-a^2}{1+\gamma}\frac{r}{1-r}+\left(\frac{1}{1+a}+\frac{r}{1-r}\right)\left(\frac{1-a^2}{1+\gamma}\right)^2\frac{r^2}{1-r^2}\\&\quad\quad-\left(\frac{1}{1+|b_0|}+\frac{r}{1+r}\right)\sum_{n=1}^{\infty}|b_n|^2r^{2n}\nonumber.
\end{align}	
	Since 
	\begin{equation}\label{e-4.4}
		\sum_{n=1}^{\infty}|a_n|r^n=\sum_{n=0}^{\infty}|b_n|r^{n+1}=r\sum_{n=0}^{\infty}|b_n|r^n,
	\end{equation}
in view of \eqref{e-4.3} and \eqref{e-4.4}, we obtain 
\begin{align*}
	\sum_{n=1}^{\infty}|a_n|r^n&\leq r\left(a+\frac{1-a^2}{1+\gamma}\frac{r}{1-r}\right)+\left(\frac{1}{1+a}+\frac{r}{1-r}\right)\left(\frac{1-a^2}{1+\gamma}\right)^2\frac{r^3}{1-r^2}\\&\quad\quad-\left(\frac{1}{1+a}+\frac{r}{1-r}\right)\sum_{n=1}^{\infty}|a_{n+1}|^2r^{2n+1}\\&=ra+\left(\frac{1-a^2}{1+\gamma}\right)\frac{r^2}{1-r}+\left(\frac{1}{1+a}+\frac{r}{1-r}\right)\left(\frac{1-a^2}{1+\gamma}\right)^2\frac{r^3}{1-r^2}\\&\quad\quad-\left(\frac{1}{1+a}+\frac{r}{1-r}\right)\sum_{n=2}^{\infty}|a_{n}|^2r^{2n-1}.
\end{align*}
Further simplification shows that
\begin{align*} &
	\sum_{n=1}^{\infty}|a_n|r^n+\left(\frac{1}{1+a}+\frac{r}{1-r}\right)\sum_{n=2}^{\infty}|a_n|^2r^{2n-1}\\&\leq ra +\left(\frac{1-a^2}{1+\gamma}\right)\frac{r^2}{1-r}+\left(\frac{1}{1+a}+\frac{r}{1-r}\right)\left(\frac{1-a^2}{1+\gamma}\right)^2\frac{r^3}{1-r^2}\\&:= \mathcal{T}(a).
\end{align*}
It is easy to see that $ \mathcal{T} $ can be represented as 
\begin{align*}
	\mathcal{T}(a)=ar+A(1-a^2)+B(1-a)(1-a^2)+C(1-a^2)^2,
\end{align*}
where 
\begin{align*}
	A=A(r)&=\frac{r^2}{(1+\gamma)(1-r)},\\ B=B(r)&=\frac{r^3}{(1+\gamma)^2(1-r^2)}\quad\mbox{and}\\ C=C(r)&=\frac{r^4}{(1+\gamma)(1-r)(1-r^2)}.
\end{align*}
Clearly, $ B $ and $ C $ are positive. We note that,
\begin{align*}
	\mathcal{T}^{\prime}(a)&=r-2Aa+B(3a^2-2a-1)+4C(a^3-a),\\\mathcal{T}^{\prime\prime}(a)&=-2A+2B(3a-1)+4C(3a^2-1) \;\; \text{and}\;\;\\ \mathcal{T}^{\prime\prime\prime}(a)&= 6B+24Ca.
\end{align*}
Since $ B $ and $ C $ are positive, it follows that $ \mathcal{T}^{\prime\prime\prime}(a)>0 $ for $ a\in [0,1] $. In other words, $ \mathcal{T}^{\prime\prime} $ is an increasing function of $ a $ in $ [0,1]. $ Therefore,
\begin{align*}
	\mathcal{T}^{\prime\prime}(a)\leq \mathcal{T}^{\prime\prime}(1)=-2A+4B+8C=\frac{2r^2}{(1+\gamma)^2(1-r)(1-r^2)}L(r),
\end{align*}
where 
\begin{align*}
	L(r)=4r^2+2r(1-r)-(1+\gamma)(1-r^2)=(1+r)(r(3+\gamma)-(1+\gamma)).
\end{align*}
It is easy to see that $ L(r)\leq 0 $ for $ r\leq r_0=(1+\gamma)/(3+\gamma) $. Hence, $ \mathcal{T}^{\prime\prime}(a)\leq 0 $ for $ a\in [0,1] $ which implies that $ \mathcal{T}^{\prime} $ is decreasing in $ [0,1] $. Therefore, for $ r\leq r_0=(1+\gamma)/(3+\gamma) $, we obtain 
\begin{align*}
	\mathcal{T}^{\prime}(a)>\mathcal{T}^{\prime}(1)=1-2A=r\frac{1+\gamma-r(3+\gamma)}{(1+\gamma)(1-r)}.
\end{align*}
Clearly, for $ r\leq r_0 $, we have $ \mathcal{T}^{\prime}(1)\geq 0 $ for all $ a\in [0,1] $. Since $ \mathcal{T}^{\prime}(a)\geq 0 $ in $ [0,1] $, $ \mathcal{T} $ is an increasing function in $ [0,1] $, and hence, we have $ \mathcal{T}(a)\leq \mathcal{T}(1)=r. $
A simple computation shows that 
\begin{align*}
	\sum_{n=0}^{\infty}|a_{n+1}|r^n+\left(\frac{1}{1+|a_1|}+\frac{r}{1-r}\right)\sum_{n=2}^{\infty}|a_n|^2r^{2(n-1)}\leq 1\;\; \text{for}\;\; r\leq r_0=\frac{1+\gamma}{3+\gamma}.
\end{align*}
To show that the sharpness of the radius we consider the function $ f_a $ by 
\begin{align*}
	f_a(z)=z\left(\frac{a-\gamma-(1-\gamma)z}{1-a\gamma-a(1-\gamma)z}\right)=B_0z-\sum_{n=1}^{\infty}B_nz^{n+1}\quad\mbox{for}\quad z\in\mathbb{D},
\end{align*}
where 
\begin{align*}
	B_0=\frac{a-\gamma}{1-a\gamma}\;\; \text{and}\;\; B_n=\frac{(1-a^2)}{a(1-a\gamma)}\left(\frac{a(1-\gamma)}{1-a\gamma}\right)^n.
\end{align*}
It is easy to see that $ a_1(f_a)=B_0 $ and $ a_n(f_a)=-B_{n-1} $. For $ n\geq 2 $, $ \gamma\in [0,1] $, and $ a>\gamma $, a simple calculation shows that 
\begin{align*}
	D(r):&=\sum_{n=1}^{\infty}|a_n|r^n+\left(\frac{1}{1+|a_1|}+\frac{r}{1-r}\right)\sum_{n=2}^{\infty}|a_n|^2r^{2n-1}\\&= \left(\frac{a-\gamma}{1-a\gamma}\right)r+\sum_{n=2}^{\infty}\frac{(1-a^2)}{a(1-a\gamma)}\left(\frac{a(1-\gamma)}{1-a\gamma}\right)^{n-1}r^n\\&\quad\quad+\left(\frac{1}{1+|B_0|}+\frac{r}{1-r}\right)\sum_{n=2}^{\infty}\frac{(1-a^2)^2}{a^2(1-a\gamma)^2}\left(\frac{a(1-\gamma)}{1-a\gamma}\right)^{2(n-1)}r^{2n-1}\\&=\left(1-\frac{1-a}{1-a\gamma}\chi(r)\right)r,
\end{align*}
where
\begin{align*}
	\chi(r)&:=1+r-\frac{(1+a)(1-\gamma)r}{1-a\gamma-a(1-\gamma)r}\\&\quad\quad-\left(\frac{1-a\gamma}{(1+a)(1-\gamma)}+\frac{r}{1-r}\right)\frac{(1+a)(1-a^2)}{1-a\gamma}\frac{(1-\gamma)^2r^2}{(1-a\gamma)^2-a^2(1-\gamma)^2r^2}.
\end{align*}
It is not difficult to show that $ \chi $ is strictly decreasing function in $ r\in (0,1) $. Hence, for $ r>r_0 $, we have $ \chi(r)<\chi(x_0) $. It is worth to point out that 
\begin{align*}
	\lim_{a\rightarrow 1}\chi(r_0)&=1+\gamma-\frac{2(1-\gamma)r_0}{1-\gamma-(1-\gamma)r_0}=1+\gamma-\frac{2r_0}{1-r_0}=0.
\end{align*}
This shows that $ \chi(r)\leq 0 $ for $ r>r_0 $ as $ a\rightarrow 1 $, and hence $ D(r)>r $ for $ r>r_0 $. Therefore
\begin{align*}
		\sum_{n=1}^{\infty}|a_n|r^n+\left(\frac{1}{1+a}+\frac{r}{1-r}\right)\sum_{n=2}^{\infty}|a_n|^2r^{2n-1}>1
\end{align*}
and hence $ r_0 $ is the best possible. This completes the proof.
\end{pf}

\noindent\textbf{Acknowledgment:}  The first author is supported by the Institute Post Doctoral Fellowship of IIT Bhubaneswar, India, the second author is supported by SERB-MATRICS, and third author is supported by CSIR, India.

\end{document}